\documentclass[a4paper,12pt]{amsart}
\usepackage{tabularx}
\usepackage[dvips]{graphicx}
\usepackage{amsmath}
\usepackage{amssymb}
\usepackage{amsbsy}
\usepackage{amsgen}
\usepackage{amsopn}
\usepackage{amscd}
\usepackage{amstext}
\usepackage{amsthm}

\theoremstyle{definition}
\newtheorem{dfn}{Definition}[section]

\theoremstyle{plain}
\newtheorem{thm}{Theorem}[section]

\theoremstyle{remark}
\newtheorem{rem}{Remark}[section]

\theoremstyle{definition}

\newtheorem{acknow}{Acknowledgment}

\newtheorem{conv}{Conventions}

\begin{document}
\title{Finite-type invariants for curves on surfaces}
\author{Noboru Ito}
\maketitle
\begin{abstract}
This study defines finite-type invariants for curves on surfaces and reveals the construction of these finite-type invariants for stable homeomorphism classes of curves on compact oriented surfaces without boundaries.  
These invariants are a higher-order generalisation of a part of Arnold's invariants that are first-order invariants for plane immersed curves.  The invariants in this theory are developed using the word theory proposed by Turaev.  
\end{abstract}

\section{Introduction.  }
V.~A.~Vassiliev developed a method for knot classification by applying the singularity theory to knots.  
This method attempts to classify knots by using  finite-type invariants (Vassiliev invariants) that are elements of $H^{0}(\mathcal{K}\setminus\Sigma)$ for the functional space $\mathcal{K}$, which comprises all the smooth mappings of $S^{1}$ into ${\mathbf{R}}^{3}$,  including the set $\Sigma$ of those mappings that are not embeddings.  
It remains unknown whether finite-type invariants can be used for the complete classification of knots (the Vassiliev conjecture).  
V.~I.~Arnold developed invariants for generic plane curves using a theory similar to that used by Vassiliev.  Arnold constructed first-order invariants that are elements of $H^{0}(\mathcal{F}\setminus\Sigma)$ for the functional space $\mathcal{F}$, which comprises all smooth mappings of $S^{1}$ into ${\mathbf R}^{2}$, including the set $\Sigma$ of mappings that are not non-generic immersions.  M.~Polyak and O.~Viro presented an explicit formula for second- or third-order Vassiliev invariants by utilizing the Gauss diagram.  Polyak also reconstructed Arnold's invariants by using the Gauss diagram in a similar manner and combinatorially defined the finite-type invariants of plane curves.  
V.~Turaev suggested that words and their topology can be considered as generalized objects of curves or knots and demonstrated that it is possible to classify knots in the same manner as words.  

In this study, the author constructs a family of finite-type invariants  for stable homeomorphism classes of curves on compact oriented surfaces without boundaries.  These invariants are a generalisation of a part of Arnold's first-order invariants and are developed using the word theory proposed by Turaev.  
In this basic paper, 
the author presents the main results and the concept of the proofs.  
The details and generalisation of the results of this study along with the background material utilized in this study will be presented in another paper.  

\section{Curves.  }  
A {\itshape curve} is a smooth immersion of an oriented circle into a closed oriented surface.  
The author now defines some types of curves.  
First, a curve is said to be {\itshape generic} if it has only transversal double points of self-intersection.  
A curve is said to be {\itshape singular} 
if it has only transversal double points, self-tangency points, and triple points of self-intersection.  
A pointed curve is defined as a curve with a base point marked along the curve except on the self-intersections.  
Two curves are said to be stably homeomorphic if their regular neighbourhoods are homeomorphic in the ambient surfaces that map the first curve onto the second one without changing the orientation of the curves and the surfaces.  
Similarly, two pointed curves are termed as stably homeomorphic if the locations of the marked base points are preserved after the curves undergo stable homeomorphism.  

\section{Definition of finite-type invariants for immersed curves.  }
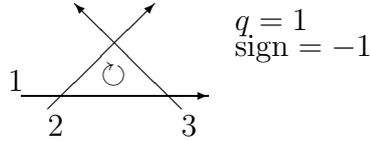
\begin{figure}
\begin{center}
\begin{picture}(70,40)
\put(0,0){\vector(1,1){40}}
\put(0,-10){$2$}
\put(-10,5){\vector(1,0){70}}
\put(-15,7){$1$}
\put(50,0){\vector(-1,1){40}}
\put(50,-10){$3$}
\put(20,10){$\circlearrowright$}
\put(70,30){$q=1$}
\put(70,20){$\mathrm{sign} = -1$}
\end{picture}
\end{center}
\caption{Sign of triangle.  }
\label{tria}
\end{figure}

Self-tangency points or triple points can be termed as {\itshape singular points}.  In particular, a self-tangency point is known as a {\itshape direct self-tangency} point if the two tangent branches are obtained using the same tangent vector; otherwise, it is called an {inverse self-tangency} point.  
The direction of the resolution of the self-tangency point is said to be {\itshape positive} if the resolution generates the part with a larger number of double points.  
The direction of the resolution of the triple point is said to be {\itshape positive} if the resolution generates a part with a {\itshape positive  triangle}, defined as the number of sides whose orientations coincide with the orientation of the curve by $q$, where the orientation of the triangle determined by a cyclic order of sides derived from three preimages of the triple point.  The sign of the triangle is defined as $(-1)^q$ (Fig.~\ref{tria}, cf.~\cite{arnold1}).  
The direction of the resolution of the singular point is said to be negative if the direction is non-positive.  
In the case of singular curves, it is possible to resolve singular points away from a marked base point.  

Let $\varphi$ denote every invariant of the surface isotopy classes of generic curves.  
The value of $\varphi$ at a singular point
is defined as the difference of the values on the curves between the positive and negative resolutions of the singular point whenever three curves are the same, except in the neighbourhood of a point where
they are as shown in each of (\ref{skein1}), (\ref{skein2}) and (\ref{skein3}): 
 
\begin{equation}
\varphi \left(
    \begin{minipage}{15pt}
        \begin{picture}(15,15)
            \qbezier(0,0)(15,7.5)(0,15)
            \qbezier(15,0)(0,7.5)(15,15)
            \put(15,15){\vector(3,2){0}}
            \put(0,15){\vector(-3,2){0}}
        \end{picture}
    \end{minipage}
\right) 
= \varphi \left(
    \begin{minipage}{15pt}
        \begin{picture}(15,15)
            \qbezier(0,0)(20,7.5)(0,15)
            \qbezier(15,0)(-5,7.5)(15,15)
            \put(15,15){\vector(2,1){0}}
            \put(0,15){\vector(-2,1){0}}
        \end{picture}
    \end{minipage}
\right)  - \varphi \left(
    \begin{minipage}{15pt}
        \begin{picture}(15,15)
            \qbezier(0,0)(10,7.5)(0,15)
            \qbezier(15,0)(5,7.5)(15,15)
            \put(15,15){\vector(1,1){0}}
            \put(0,15){\vector(-1,1){0}}
        \end{picture}
    \end{minipage}
\right), \label{skein1}
\end{equation}
\begin{equation}
\varphi \left(
    \begin{minipage}{15pt}
        \begin{picture}(15,15)
            \qbezier(0,0)(15,7.5)(0,15)
            \qbezier(15,0)(0,7.5)(15,15)
            \put(15,15){\vector(3,2){0}}
            \put(0,0){\vector(-3,-2){0}}
        \end{picture}
    \end{minipage}
\right) 
= \varphi \left(
    \begin{minipage}{15pt}
        \begin{picture}(15,15)
            \qbezier(0,0)(20,7.5)(0,15)
            \qbezier(15,0)(-5,7.5)(15,15)
            \put(15,15){\vector(2,1){0}}
            \put(0,0){\vector(-2,-1){0}}
        \end{picture}
    \end{minipage}
\right)  - \varphi \left(
    \begin{minipage}{15pt}
        \begin{picture}(15,15)
            \qbezier(0,0)(10,7.5)(0,15)
            \qbezier(15,0)(5,7.5)(15,15)
            \put(15,15){\vector(1,1){0}}
            \put(0,0){\vector(-1,-1){0}}
        \end{picture}
    \end{minipage}
\right), \label{skein2}
\end{equation}
\begin{equation}
\varphi \left(
    \begin{minipage}{15pt}
        \begin{picture}(15,15)
            \put(0,0){\line(1,1){15}}
           \put(15,0){\line(-1,1){15}}
            \put(7.5,-1){\line(0,1){18}}
            \end{picture}
    \end{minipage}
\right)  
= \varphi \left(
    \begin{minipage}{15pt}
        \begin{picture}(15,15)
            \put(0,0){\line(1,1){15}}
           \put(15,0){\line(-1,1){15}}
            \put(3,-1){\line(0,1){18}}
            \end{picture}
    \end{minipage}
\right) - \varphi \left(
    \begin{minipage}{15pt}
        \begin{picture}(15,15)
            \put(0,0){\line(1,1){15}}
           \put(15,0){\line(-1,1){15}}
            \put(12,-1){\line(0,1){18}}
            \end{picture}
    \end{minipage}
\right).  \label{skein3}
\end{equation}

\begin{minipage}{15pt}
\begin{picture}(100,5)
\put(180,5){positive} \put(230,5){negative}
\end{picture}
\end{minipage}

Every invariant of the curves is extended inductively to that of singular curves by resolving the singular points using (\ref{skein1}), (\ref{skein2}), and (\ref{skein3}).  

\begin{dfn}
An invariant of generic curves $\varphi$ is said to be 
{\itshape a finite-type invariant of order less than or equal to $n$} if $\varphi$ vanishes on every singular curve with at least $n+1$ singular points (self-tangency points or triple points), where $\varphi$ is expanded by (\ref{skein1}), (\ref{skein2}) and (\ref{skein3}).  
\end{dfn}

\section{Signed words.  }
\begin{dfn}
Let $\tau$ be the involution: $2j$ $\mapsto$ $2j-1$ and $2j-1$ $\mapsto$ $2j$ for every $j \in {\mathbf N}$.  
Let $\mathcal{A}$ be the set comprising $\{$$X_{j}, \overline{X}_{j}$ $($$j \in \mathbf{N}$$)$ $|$ $X_{j}$ $:=$ $\tau(2j)$, $\overline{X}_{j}$ $:=$ $2j$, $\overline{\overline{X}}_{j} = X_{j}$$\}$.  A {\itshape signed word of length $2n$} is a mapping from $\hat{n}:= \{1, 2, 3, \dots, 2n\}$ to $\mathcal{A}$ where each element of $w(\hat{n})$ is the image of precisely two elements of $\hat{n}$ and $w(\hat{n})$ has not both $X_{j}$ and $\overline{X}_{j}$ for every $j \in {\mathbf N}$.  
A signed word of length $0$ is denoted by $\emptyset$.  For every signed word $w$ of length $2n$, the elements of $w(\hat{n})$ are termed as {\itshape letters} of the signed word $w$.  
For every letter $A$ of a signed word $w$, set ${{\mathrm sign}_{w}} A$ $=$ $-1$ if $A = \overline{X}$ and ${\mathrm sign}_{w} A$ $=$ $1$ if $A =X$.  
Two signed words $w$ and $w'$ of length $2n$ are {\itshape isomorphic} 
if there is a bijection $f$ $:$ ${\mathcal A} \to {\mathcal A}$ such that $w'$ $=$ $f w$ and ${{\mathrm sign}_{w'}} w'(i)$ $=$ ${{\mathrm sign}_{f w}} f w(i)$ for every $i \in \hat{n}$.  The isomorphism of two signed words $w$ and $w'$ is denoted by $w$ $\simeq$ $w'$.  For example, $\overline{X}_{1}X_{2}X_{2}\overline{X}_{1}$ $\simeq$ $\overline{X}_{5}X_{2}X_{2}\overline{X}_{5}$.  However, $\overline{X}_{1}X_{2}X_{2}\overline{X}_{1}$ is not isomorphic to $X_{1}\overline{X}_{2}\overline{X}_{2}X_{1}$.  
For two signed words $u$ and $w$, $u$ $\prec$ $w$ implies that $u$ is a subword of $w$.  
\end{dfn}

\section{Construction of the invariant.  }\label{const-in}
For two arbitrary signed words $u$ and $w$, 
$\langle \ , \rangle$ can be defined by 
\begin{equation}\label{def-bilin}
\left\langle u, w \right\rangle = \sum_{v \prec w} (u, v), 
\end{equation}
where $(u, v)$ is $1$ if $u$ $\simeq$ $v$ and is $0$ otherwise.  
Let $k$ be a field; ${\mathbf W}$, the $k$-linear space generated by all the isomorphic classes of the signed words; and ${\mathbf W}^{*}$, the dual space of ${\mathbf W}$.  
We expand $\langle \ , \rangle$ linearly to $\langle \ , \rangle$ $:$ ${\mathbf W} \times {\mathbf W} \to k$.  

For an arbitrary generic curve $\Gamma$, which has  $m$ double points, 
let $w_\Gamma$ be a signed word$: \hat{m} \to \mathcal{A}$ that is determined by selecting an arbitrary marked base point in the following manner as in \cite{turaev2}.  

For a given generic curve $\Gamma$, which has $m$ double points, used to define $w_{\Gamma}$, first,  the double points of $\Gamma$ are labelled by using distinct letters $A_{1}, A_{2}, \dots, A_{2m}$ such that $A_{i}$ $=$ $\overline{X}_{i}$ or $X_{i}$.  
Then, beginning with a marked base point and moving along $\Gamma$ until returning to the base point, the first instances of all the double points are labelled.  Let the $i$-th double point be labelled $A_{i}$.  Let $t^{1}_{i}$ (resp. $t^{2}_{i}$) be the tangent vector to $\Gamma$ at the first (resp. second) passage through this double point.  Set $A_{i} = \overline{X}_{i}$ if the pair $(t^{1}_{i}, t^{2}_{i})$ is \begin{minipage}{30pt}
        \begin{picture}(30,30)
            \put(5,10){\vector(1,1){20}}\put(-5,2){\small $1 {\rm st}$}
           \put(25,10){\vector(-1,1){20}}\put(20,2){\small $2 {\rm nd}$}
            \end{picture}
             \end{minipage} and $A_{i} = X_{i}$ if otherwise.  Since every double point is traversed twice, this yields a signed word $w_{\Gamma}$, which is represented by $m$ elements in $\mathcal{A}$.

\begin{dfn}\label{alter}
The $\nu$-shift is defined as $\nu(AxAy)$ $=$ $x\overline{A}y\overline{A}$, where 
$A \in \mathcal{A}$ and $x$ and $y$ are subwords.  
The {\itshape cyclic equivalence} $\sim$ is defined as 
$w \sim w'$ if and only if there exists $l \in {\mathbf N}$ such that $\nu^{l}(w) = w'$.  
Let ${\mathbf W}_{n}$ be the $k$-linear space generated by all the isomorphic classes with signed words of length $2n$ and $\mathbf W_{n}^{*}$ be the dual space of $\mathbf W_{n}$.  Let the subspace ${\mathbf W}_{n}^{\nu}$ of ${\mathbf W}_{n}$ be the vector space generated by $\{w \in {\mathbf W}_{n }|~\nu(w) = w \}$.  Let $({\mathbf W}_{n}^{\nu})^{*}$ be the dual space of ${\mathbf W}_{n}^{\nu}$.  
For every cyclic equivalence class containing a signed word $v$, the sum of all the representative elements of the cyclic equivalence class is denoted by $[v]$.  For two arbitrary $\alpha$ and $\beta$ $\in k$ and two arbitrary $u$ and $v$ $\in {\mathbf W}$, let $[\alpha u + \beta v]$ $:=$ $\alpha[u]$ $+$ $\beta[v]$.  
\end{dfn}

For example, if the cyclic equivalence of the isomorphic class of the signed words is $\{$$X_{1}X_{1}\overline{X}_{2}\overline{X}_{2}$, $\overline{X}_{1}\overline{X}_{2}\overline{X}_{2}\overline{X}_{1}$, $\overline{X}_{2}\overline{X}_{2}X_{1}X_{1}$, $X_{2}X_{1}X_{1}X_{2}$$\}$ and $v$ $=$ $\overline{X}_{2}\overline{X}_{2}X_{1}X_{1}$, then
$[v]$ $=$ $X_{1}X_{1}\overline{X}_{2}\overline{X}_{2}$ $+$ $\overline{X}_{1}\overline{X}_{2}\overline{X}_{2}\overline{X}_{1}$ $+$ $\overline{X}_{2}\overline{X}_{2}X_{1}X_{1}
$ $+$ $X_{2}X_{1}X_{1}{X}_{2}$.  

We denote the linear space generated by all stable homeomorphism classes of curves on compact oriented surfaces without boundaries by $\mathcal{C}$.  
\begin{rem} 
By using $\overline{X}$ $\mapsto$
             \begin{minipage}{30pt}
        \begin{picture}(30,30)
            \put(0,10){\vector(1,1){20}}\put(-5,3){\small $1 {\rm st}$}
           \put(20,10){\vector(-1,1){20}}\put(20,3){\small $2 {\rm nd}$, }
            \end{picture}
             \end{minipage}
\quad and \quad $X$ $\mapsto$
\begin{minipage}{30pt}
        \begin{picture}(30,30)
            \put(0,10){\vector(1,1){20}}\put(-5,3){\small $2 {\rm nd}$}
           \put(20,10){\vector(-1,1){20}}\put(20,3){\small $1 {\rm st}$}
            \end{picture}
             \end{minipage}\quad
every signed word determines a regular neighbourhood of a curve $\Gamma$ on a surface ${\mathcal S}$, where $\Gamma$ gives the CW-decomposition of ${\mathcal S}$ (Fig.~\ref{reg-nbd}).  
We denote the linear space generated by all stable homeomorphism classes of curves with $n$ double points on compact oriented surfaces without boundaries by $\mathcal{C}_{n}$.  
There exists a bijective mapping from $\mathcal{C}_{n}$ to ${\mathbf W}_{n}^{\nu}$; this has been proved by V. Turaev \cite{turaev2}.  
In the rest of this paper, we identify $\mathcal{C}_{n}$ with ${\mathbf W}_{n}^{\nu}$ and $\mathcal{C}$ with ${\mathbf W}^{\nu}$, where ${\mathbf W}^{\nu}$ is the vector space generated by $\{w \in {\mathbf W}|~\nu(w) = w \}$.  
\begin{figure}
\[\]
\[\overline{X}_{1}X_{2}\overline{X}_{1}X_{2} \leftrightarrow
\begin{minipage}{80pt}
\begin{picture}(50,50)
\put(55,37){$X_{2}$}
\put(57,16){$\overline{X}_{1}$}
\includegraphics[width=3cm,height=3cm]{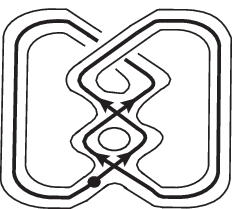}
\end{picture}
\end{minipage}
\quad
\cup {\text{ two disks}}.  \]
\caption{Signed words and pointed curves.  }\label{reg-nbd}
\end{figure}
\end{rem}

For a given generic curve $\Gamma$, let $w_\Gamma$ be a signed word that is determined by  $\Gamma$.  For two arbitrary natural numbers $m$ and $n$, 
we define a {\it signed curve invariant of order} $n$, $SCI_{n}$ $:$ ${\mathbf W}_{m}^{\nu}$ $\to$ $({\mathbf W}_{n}^{\nu})^{*}$, as follows: 
\begin{equation}\label{def-sci}
SCI_n(\Gamma) = \left\langle [\cdot] , w_{\Gamma} \right\rangle: {\mathbf W}_{n} \to k.  
\end{equation}
By the definition of $SCI_{n}$, we use the following notation: 
\begin{equation}\label{rep-sci}
SCI_{n}(\Gamma)(v) = \left\langle [v] , w_{\Gamma} \right\rangle \ \ \ (v \in {\mathbf W}_{n}).  
\end{equation}
By using the definition of $[~]$ and $\langle \ , \rangle$,  $SCI_{n}$ is independent of the choice of the marked base point on the curve $\Gamma$.  

\begin{thm}\label{main}
$SCI_n$ is a finite-type invariant of order less than or equal to $n$ for an arbitrary generic curve on a surface.  
\end{thm}
The proof of theorem \ref{main} is given in Sect.~\ref{main}.  

Let $m_{w}$ be the number of representative elements of the cyclic equivalence class that contains a signed word $w$.  Let $\overline{w}$ $=$ $\frac{1}{m_{w}} [w]$ $\in {\mathbf W}^{\nu}$.  Then, 
$SCI_{n}(\Gamma)$ $=$ $\langle [\cdot], \overline{w}_{\Gamma} \rangle$.  
For two arbitrary signed words $v$ and $u$, 
$v^{*}$ denotes a linear mapping such that $v^{*}(u)$ $=$ $(v, u)$.  
$SCI_{n}$ restricted on ${\mathbf W}_{n}^{\nu}$ gives  the isomorphism between ${\mathbf W}_{n}^{\nu}$ and $({\mathbf W}_{n}^{\nu})^{*}$, sending $\overline{w}_{\Gamma}$ to $[w_{\Gamma}]^{*}$.  
We denote this isomorphism by $\iota_{n}$.  
\begin{thm}
For $SCI_{k}$ $:$ ${\mathbf W}^{\nu}_{n} \to ({\mathbf W}^{\nu}_{k})^{*}$, 
we establish the following relation: 
\begin{equation}\label{sci-eq}
\begin{split}
\frac{(n-k + 1)!}{(n-l-1)!(k-l)!}SCI_{l} = SCI_{l}\circ&\iota_{k}^{-1} \circ SCI_{k}\\ &(1\le l \le k \le n).  
\end{split}
\end{equation}
\end{thm}
\begin{equation}\label{sci-eq2}
(n-k + 1)SCI_{k-1} = SCI_{k-1}\circ \iota_{k}^{-1} \circ SCI_{k} \ (2\le k \le n) 
\end{equation}
leads to (\ref{sci-eq}).  

\section{Arnold's invariants and $SCI_{n}$.  }\label{ar-sci}
In Sect.~\ref{ar-sci}, $\Gamma$ is used to denote an arbitrary plane generic curve.  Let $w_{\Gamma}$ be a signed word determined by selecting an arbitrary marked base point by using the definition in Sect.~\ref{const-in}.  
Arnold's basic invariants $(J^+, J^-, St)$ and $SCI_{n}$ ($n \in \mathbf{N}$) are invariants for plane isotopy classes of curves on a plane.  
The definitions of $J^+$, $J^-$ and $St$ are provided in \cite{arnold2, arnold1}.  

Let $i(\Gamma)$ be the rotation number for every plane generic curve $\Gamma$ and then, $i(\Gamma)$ is a finite-type invariant of order $0$ because the left-hand sides of (\ref{skein1}), (\ref{skein2}) and (\ref{skein3}) should be equal to $0$.  
We obtain the following equations: 
\begin{equation}
J^{+}(\Gamma) - J^{-}(\Gamma) = SCI_{1}(\Gamma)(X_{1}X_{1}), 
\end{equation}
\begin{equation}
\begin{split}
J^{-}(\Gamma) + 6St(\Gamma) = &- 2 SCI_{2}(\Gamma)(X_{1}X_{1}X_{2}X_{2}\\&
- X_{1}X_{1}\overline{X}_{2}\overline{X}_{2} + \overline{X}_{1}\overline{X}_{1}\overline{X}_{2}\overline{X}_{2})\\
&+ i^{2}(\Gamma) -1.  
\end{split}
\end{equation}

However, $J^{+}$, $J^{-}$ and $St$ cannot be represented by $i(\Gamma)$, $SCI_{1}$ and $SCI_{2}$ in a similar manner.  

\section{Proof of the main theorem.  \label{main}}
It is now shown that for an expanded invariant $I$ for singular curves with $SCI_n$, $I(\Gamma)$ $=$ $0$ for an arbitrary singular curve $\Gamma$ with $m$ ($\ge n+1$) singular points ($P_1, P_2, \dots, P_m$).  
We define $\sigma$ by the mapping  $P_i$ $\mapsto$ $\sigma_i$ and $\sigma$ $=$ $($$\sigma_1, \sigma_2, \dots, \sigma_m$$)$, where $\sigma_{i} = \pm 1$.  Let $\mathcal{J}$ $=$ $\{1, 2, \dots, m\}$ and $\mathrm{sign}\,\sigma$ $=$ $\prod_{i \in \mathcal{J}}\sigma_i$.  
The generic curve is obtained from  $\Gamma$ by resolving at $P_i$, where $i = 1, 2, \dots, m$, in the positive direction if $\sigma_i$ $=$ $1$ and in the negative direction if $\sigma_i$ $=$ $-1$, and it is denoted by $\Gamma_{\sigma}$.  Let $w_{\Gamma_{\sigma}}$ be a signed word derived from $\Gamma_{\sigma}$ by the selection of an arbitrary marked base point.  

By using the definitions of $\sigma$ and $\Gamma_{\sigma}$ in (\ref{skein1}), (\ref{skein2}) and (\ref{skein3}), we obtain 
\begin{equation}\label{sci-sin}
I(\Gamma) = \displaystyle \sum_{\sigma \in \{-1, 1\}^{m}}\mathrm{sign}\,\sigma \,SCI_n(\Gamma_{\sigma}).  
\end{equation}
Then, by using (\ref{def-bilin}), (\ref{rep-sci}), (\ref{sci-sin}) and $[v]$ $=$ $\sum_{u}$ $u$, where $v$ $\sim$ $u$, we obtain 
\begin{equation}\label{coeffi}
\begin{split}
I(\Gamma)(v)
&= \sum_{\sigma \in \{-1, 1\}^{m}} \mathrm{sign}\,\sigma \sum_{u (\sim v)} \left \langle u, w_{\Gamma_{\sigma}} \right \rangle \\
&= \sum_{u (\sim v)} \sum_{\sigma \in \{-1, 1\}^{m}}\mathrm{sign}\,\sigma \left \langle u, w_{\Gamma_{\sigma}} \right \rangle.  
\end{split}
\end{equation}

In order to show that $\sum_{\sigma}$ $\mathrm{sign}\,\sigma$ $\left \langle u, w_{\Gamma_{\sigma}} \right \rangle$ vanishes, we consider the following conditions.  
For an arbitrary $\Gamma$, we can determine the set of generic curves $\{ \Gamma_{\sigma} \}_{\sigma}$.  
Let $w_{\sigma_{0}}$ be a signed word corresponding to $\Gamma_{\sigma_{0}}$ with $\sigma_{0}$ $=$ $(1, 1, \dots, 1)$.  By using $w_{\sigma_{0}}$, $w_{\sigma}$ is determined as follows.  

Let $A$, $B$, $C$, $\overline{A}$, $\overline{B}$ and $\overline{C}$ be the letters of signed words, and $x$, $y$, $z$ and $t$, letters of subwords.  
In the case of a negative resolution of a direct self-tangency point, it is necessary to consider two cases: 
if $w_{\Gamma_{\sigma_{0}}}$ $=$ $x\overline{A}By\overline{A}Bz$, then $w_{\Gamma_{\sigma}}$ $=$ $xyz$ (Fig.~\ref{direct-self}), and if $w_{\Gamma_{\sigma_{0}}}$ $=$ $xA\overline{B}yA\overline{B}z$, then $w_{\Gamma_{\sigma}}$ $=$ $xyz$.  
Even in the case of a negative resolution of an  inverse self-tangency point, there are two cases: 
if $w_{\Gamma_{\sigma_{0}}}$ $=$ $x\overline{A}ByB\overline{A}z$, then $w_{\Gamma_{\sigma}}$ $=$ $xyz$, and if $w_{\Gamma_{\sigma_{0}}}$ $=$ $xA\overline{B}y\overline{B}Az$, then $w_{\Gamma_{\sigma}}$ $=$ $xyz$ (Fig.~\ref{inverse-self}).  
In the case of a negative resolution of a triple point, there are eight cases: 
if $w_{\Gamma_{\sigma_{0}}}$ $=$ $xAByCBzCAt$, then $w_{\Gamma_{\sigma}}$ $=$ $xBAyBCzACt$; 
if $w_{\Gamma_{\sigma_{0}}}$ $=$ $x\overline{A}By\overline{A}CzCBt$, then $w_{\Gamma_{\sigma}}$ $=$ $xB\overline{A}yC\overline{A}zBCt$; 
if $w_{\Gamma_{\sigma_{0}}}$ $=$ $xA\overline{B}yCAz\overline{B}Ct$, $w_{\Gamma_{\sigma}}$ $=$ $x\overline{B}AyACzC\overline{B}t$ (Fig.~\ref{triple-resolution}); 
if $w_{\Gamma_{\sigma_{0}}}$ $=$ $xAByB\overline{C}zA\overline{C}t$, then $w_{\Gamma_{\sigma}}$ $=$ $xBAy\overline{C}Bz\overline{C}At$; and four more patterns, in which each of the three letters take an additional bar to represent the different resolutions of the situations given above.  

\begin{figure}[t]
\begin{center}
    \begin{minipage}{30pt}
        \begin{picture}(40,40)
            \qbezier(0,10)(40,25)(0,40)
            \qbezier(30,10)(-10,25)(30,40)
            \put(30,40){\vector(3,1){0}}
            \put(0,40){\vector(-3,1){0}}
            \put(20,28){$B$}
            \put(22,16){$\overline{A}$}
        \end{picture}
     \end{minipage}
\quad
$\stackrel{\sigma_{i} = 1}{\longleftarrow}$
\quad
\begin{minipage}{30pt}
        \begin{picture}(40,40)
            \qbezier(0,10)(30,25)(0,40)
            \qbezier(30,10)(0,25)(30,40)
            \put(30,40){\vector(4,3){0}}
            \put(0,40){\vector(-4,3){0}}
        \end{picture}
    \end{minipage}
\quad
$\stackrel{\sigma_{i} = -1}{\longrightarrow}$
\quad
    \begin{minipage}{30pt}
        \begin{picture}(40,40)
            \qbezier(0,10)(20,25)(0,40)
            \qbezier(30,10)(10,25)(30,40)
            \put(30,40){\vector(3,2){0}}
            \put(0,40){\vector(-3,2){0}}
            \end{picture}
    \end{minipage} 
\end{center}
\  
$\Gamma_{\sigma_{0}}$\qquad\qquad\qquad$\Gamma$\quad\qquad\qquad\qquad\ $\Gamma_{\sigma}$\\
$w_{\Gamma_{\sigma_{0}}} = x\overline{A}By\overline{A}Bz$
\qquad\qquad\qquad\qquad\quad $w_{\Gamma_{\sigma}} = xyz$
\caption{An example of $w_{\Gamma_{\sigma}}$ 
obtained by a negative resolution at a direct self-tangency point $P_i$}\label{direct-self}
\end{figure}
\begin{figure}[t]
\begin{center}
    \begin{minipage}{30pt}
        \begin{picture}(40,40)
            \qbezier(0,10)(40,25)(0,40)
            \qbezier(30,10)(-10,25)(30,40)
            \put(30,40){\vector(3,1){0}}
            \put(0,10){\vector(-3,-1){0}}
            \put(23,27){$\overline{B}$}
            \put(19,16){$A$}
        \end{picture}
    \end{minipage}
\quad
$\stackrel{\sigma_{i} = 1}{\longleftarrow}$
\quad
\begin{minipage}{30pt}
        \begin{picture}(40,40)
            \qbezier(0,10)(30,25)(0,40)
            \qbezier(30,10)(0,25)(30,40)
            \put(30,40){\vector(3,2){0}}
            \put(0,10){\vector(-3,-2){0}}
        \end{picture}
    \end{minipage}
\quad
$\stackrel{\sigma_{i} = -1}{\longrightarrow}$
\quad
    \begin{minipage}{30pt}
        \begin{picture}(40,40)
            \qbezier(0,10)(20,25)(0,40)
            \qbezier(30,10)(10,25)(30,40)
            \put(30,40){\vector(1,1){0}}
            \put(0,10){\vector(-1,-1){0}}
        \end{picture}
    \end{minipage} 
\end{center}
\  
$\Gamma_{\sigma_{0}}$\qquad\qquad\qquad$\Gamma$\quad\qquad\qquad\qquad\ $\Gamma_{\sigma}$\\
$w_{\Gamma_{\sigma_{0}}} = xA\overline{B}y\overline{B}Az$
\qquad\qquad\qquad\qquad\quad $w_{\Gamma_{\sigma}} = xyz$
\caption{An example of $w_{\Gamma_{\sigma}}$ obtained by a negative resolution at an inverse self-tangency point $P_i$}\label{inverse-self}
\end{figure}
\begin{figure}[t]
\begin{center}
        \begin{minipage}{30pt}
        \begin{picture}(40,40)
            \put(0,10){\vector(1,1){30}}
           \put(30,10){\vector(-1,1){30}}
            \put(6,42){\vector(0,-1){34}}
            \put(18,23){$A$}
            \put(8,35){$\overline{B}$}
            \put(-4,13){$C$}
            \end{picture}
    \end{minipage}
\quad
$\stackrel{\sigma_{i} = 1}{\longleftarrow}$
\quad
\begin{minipage}{30pt}
        \begin{picture}(40,40)
            \put(0,10){\vector(1,1){30}}
           \put(30,10){\vector(-1,1){30}}
            \put(15,42){\vector(0,-1){32}}
            \end{picture}
    \end{minipage}
\quad
$\stackrel{\sigma_{i} = -1}{\longrightarrow}$
\quad
    \begin{minipage}{30pt}
        \begin{picture}(40,40)
            \put(0,10){\vector(1,1){30}}
           \put(30,10){\vector(-1,1){30}}
            \put(24,42){\vector(0,-1){34}}
            \put(2,23){$A$}
            \put(15,32){$C$}
            \put(27,14.8){$\overline{B}$}
            \end{picture}
    \end{minipage}\\
\  
$\Gamma_{\sigma_{0}}$\qquad\qquad\qquad$\Gamma$\quad\qquad\qquad\qquad\ $\Gamma_{\sigma}$\\
$w_{\Gamma_{\sigma_{0}}} = xA\overline{B}yCAz\overline{B}Ct$
\qquad $w_{\Gamma_{\sigma}} = x\overline{B}AyACzC\overline{B}t$  
\caption{An example of $w_{\Gamma_{\sigma}}$ obtained by a negative resolution at a triple point $P_i$}\label{triple-resolution}
\end{center}
\end{figure}

Let $u'$ denote a subword of $w_{\Gamma_{\sigma_{0}}}$ and $\sigma(u')$, a subword of $w_{\Gamma_{\sigma}}$ corresponding to $u'$.  $\sigma(u')$ can be considered to be $\emptyset$ for a negative resolution of a self-tangency point.  Let $v'$ be a signed word.  The subword $v'$ is represented as $\sigma(u')$.  
By using (\ref{def-bilin}) and the notion above, 
\begin{equation}\label{exchan}
\begin{split}
&\sum_{\sigma \in \{-1, 1\}^{m}} \mathrm{sign}\,\sigma \left \langle u, w_{\Gamma_{\sigma}} \right \rangle \\
&= \sum_{\sigma \in \{-1, 1\}^{m}} \mathrm{sign}\,\sigma \sum_{v' \prec w_{\Gamma_{\sigma}}} (u, v') \\
&= \sum_{\sigma \in \{-1, 1\}^{m}} \mathrm{sign}\,\sigma \sum_{\sigma(u') \prec w_{\Gamma_{\sigma}}} (u, \sigma(u')) \\
&= \sum_{\sigma \in \{-1, 1\}^{m}} \mathrm{sign}\,\sigma \sum_{u' \prec w_{\Gamma_{\sigma_{0}}}} (u, \sigma(u')) \\
&= \sum_{u' \prec w_{\Gamma_{\sigma_{0}}}} \sum_{\sigma \in \{-1, 1\}^{m}} \mathrm{sign}\,\sigma (u, \sigma(u')).  
\end{split}
\end{equation}
Let $\mathcal{J}_{u'}$ be $\{$ $i \in {\mathcal J}$ $|$ a  letter of $\sigma(u')$ is generated by a resolution at $P_i$ $\}$.  
The sum concerning $\sigma$ can be divided into one part with $\sigma_i$ $(i \in \mathcal{J}_{u'})$ and another part with $\sigma_i$ $(i \in \mathcal{J} \setminus \mathcal{J}_{u'})$.  
Hence, we have 
\begin{equation}\label{decompo}
\begin{split}
&\sum_{\sigma \in \{-1, 1\}^{m}} {\mathrm{sign}}\,\sigma (u, \sigma(u')) \\
&= \sum_{
\begin{subarray}{c}
\sigma_{i} = \pm 1 \\
i \in \mathcal{J}_{u'}
\end{subarray}
} \prod_{i \in \mathcal{J}_{u'}} \sigma_{i} \sum_{
\begin{subarray}{c}
\sigma_{i} = \pm 1 \\
i \in \mathcal{J} \setminus \mathcal{J}_{u'}
\end{subarray}
} \prod_{i \in \mathcal{J} \setminus \mathcal{J}_{u'}} \sigma_{i} (u, \sigma(u')).  
\end{split}
\end{equation}

For an arbitrary $u$, where $\sigma(u')$ $\in$ ${\mathbf W}_{n}$, there exists at least one singular point $P_i$ that is not related to $(u, \sigma(u'))$ because $u$ and $u'$ $\in$ ${\mathbf W}_{n}$ and $m$ $\ge n+1$.  In other words, 
$\mathcal{J} \setminus \mathcal{J}_{u'}$ is not empty.  
Then, 

\begin{equation}\label{last-eq}
\begin{split}
&\sum_{
\begin{subarray}{c}
\sigma_{i} = \pm 1 \\
i \in \mathcal{J} \setminus \mathcal{J}_{u'}
\end{subarray}
} \prod_{i \in \mathcal{J} \setminus \mathcal{J}_{u'}} \sigma_{i} (u, \sigma(u')) \\
&= (u, \sigma(u'))
\sum_{
\begin{subarray}{c}
\sigma_{i} = \pm 1 \\
i \in \mathcal{J} \setminus \mathcal{J}_{u'}
\end{subarray}
} \prod_{i \in \mathcal{J} \setminus \mathcal{J}_{u'}} \sigma_{i} \\
&= 0.  
\end{split}
\end{equation}
Here, (\ref{last-eq}) implies that $I(\Gamma)(v)$ in (\ref{coeffi}) vanishes, and $I(\Gamma)$ $=$ $0$.

\begin{acknow}
The author is grateful to Professors K. Habiro, T. Ohtsuki, Y. Ohyama and K. Taniyama for their comments on this study.  
The author also thanks M. Fujiwara for his useful comments about the relation to Arnold's invariants in an earlier version of this paper.  
Finally, the author expresses his deepest gratitude to his advisor Professor J. Murakami for his encouragement and advice.  
\end{acknow}

\end{document}